\documentclass{article}
\usepackage{amscd,amssymb,hyperref,mathtools}
\usepackage{amsmath}
\usepackage{cite}
\usepackage[utf8]{inputenc}
\usepackage[all]{xy}
\usepackage{amsmath}
\usepackage{graphicx}
\usepackage{enumerate}
\usepackage{verbatim}
\newtheorem{lemma}{Lemma}

\newcommand{\R}{\mathbb{R}}

\newcommand{\C}{\mathbb{C}}

\newcommand{\Hyp}{\mathbb{H}}

\newcommand{\be}{\begin{enumerate}}
\newcommand{\ee}{\end{enumerate}}
\newcommand{\bi}{\begin{itemize}}
\newcommand{\ei}{\end{itemize}}

\newcommand{\ba}{\begin{array}}
\newcommand{\ea}{\end{array}}
\newcommand{\bmat}{\left[\begin{array}}
\newcommand{\emat}{\end{array}\right]}
\newcommand{\bt}{\begin{thm}}
\newcommand{\et}{\end{thm}}
\newcommand{\bp}{\begin{proof}}
\newcommand{\ep}{\end{proof}}
\newcommand{\bprop}{\begin{prop}}
\newcommand{\eprop}{\end{prop}}
\newcommand{\bl}{\begin{lemma}}
\newcommand{\el}{\end{lemma}}
\newcommand{\bc}{\begin{cor}}
\newcommand{\ec}{\end{cor}}
\newcommand{\bd}{\begin{defn}}
\newcommand{\ed}{\end{defn}}
\newcommand{\ov}{\overline}

\newcommand{\vep}{\varepsilon}

\newcommand{\tx}{\textrm}

\title{Conformal Group Actions on Generalized Kuramoto Oscillators}
\author{Max Lipton  \\
	Cornell University  \\
	}

\date{December 19, 2018}
\begin{document}

\maketitle

\begin{abstract}
\noindent This paper unifies the recent results on generalized Kuramoto Model reductions. Lohe took a coupled system of $N$ bodies on $S^d$ governed by the Kuramoto equations $\dot{x_i} = \Omega x_i + X - \langle x_i, X \rangle x_i$ and used the method of Watanabe and Strogatz to reduce this system to $d + \frac{d(d-1)}{2}$ equations. Using a model of rigid rotations on a sphere as a guide, we show that the reduction is described by a smooth path in the Lie group of conformal transformations on the sphere, which is diffeomorphic to $SO(d) \times D^d$. Seeing the reduction this way allows us to apply geometric and topological reasoning in order to understand qualitative behavior of the Kuramoto Model. \\

\noindent keywords: Kuramoto Model, dynamical systems, Lie groups, conformal geometry
\end{abstract}

\normalsize

\section{Introduction}

Suppose we mark $N$ points on a sphere $S^d$ and rotate it over time via a continuous path $R(t)$ in $SO(d+1)$, tracking the location of the individual points. Noether's Theorem tells us there is a system of conserved quantities, and it is readily apparent the angles between the points (in $\R^{d+1}$ from the standard sphere embedding) are conserved over time. Rather than expressing the position of point $x_i$ based on the given rotations of the sphere, suppose we start with the $N$ coupled differential equations for the points and we deduce the pairwise angles between points, or some equivalent quantity like the pairwise inner products, are conserved. Then from this information, we can derive a single differential equation $\dot{R} = f(R)$ in $SO(d+1)$, dramatically reducing the number of variables. The Kuramoto Model follows this footprint, except the conserved quantities are the cross ratios preserved by conformal transformations.

Mirollo connected the classical model on $S^1$ to the M\"obius group of fractional linear transformations on $\R^2$, a topic which is well-understood\cite{mirollo}. Max Lohe derived reduced equations in higher dimensional Kuramoto models\cite{lohe}, whilst Chandra, Girvan, and Ott derive the higher dimensional continuum limit on the number of bodies\cite{ott}. The conserved quantities Lohe derived are those that are preserved by conformal transformations on the sphere. We have tied up these disparate results into a unified statement about conformal dynamical systems, and hopefully we can use this interpretation to get a concrete qualitative statement on higher dimensional Kuramoto Models.

\section{Classical Kuramoto Results}

We denote the $d$-dimensional sphere by $S^d$, and we denote the $d$-dimensional open disc by $D^d$.

The Kuramoto Model is a well-studied collection of dynamical systems which model coupled oscillators. Depending on the choice of coupling and its parameters, phenomena such as synchronization, segregation, and scattering can occur. Kuramoto Models have a diversity of applications including population sleep cycles, Josephson junctions, systems of mechanically linked motors, and more.

The classical one dimensional case concerns a system of $N$ coupled oscillators on $S^1 \subseteq \C$, written as $\dot{\theta_j} = A + B\cos{\theta_j} + C\sin{\theta_j}, j = 1 \dots N$, where $A, B,$ and $C$ are smooth functions of time and the angular variables $\theta_1, \dots, \theta_N$. For many of the preexisting results on conserved quantities to hold, we require that the coefficient functions are the same for all $j$ and they involve all the angular variables. Sometimes the system will be expressed in a form that obfuscates $A, B,$ and $C$, but by some clever use of trigonometric identities, we can rewrite the governing equations in the desired form. For example, the first moment system $\dot{\theta_j} = \omega + \frac{1}{N} \sum\limits_{i=1}^N e^{i(\theta_i - \theta_j)}$ ostensibly has different coefficient functions which may vary with $j$, but observe we may rewrite $\dot{\theta_j} = \omega + \Big(\frac{1}{N}\sum\limits_{i = 1}^N e^{i\theta_i} \Big) \cos{\theta_j} + \Big(\frac{1}{iN}\sum\limits_{i = 1}^N e^{i \theta_i} \Big) \sin{\theta_j}$. If the functions are invariant under permutations of the $\theta_j$, we say the system is symmetric, like the aforementioned example.

Lohe describes how to incorporate different amplitudes, frequencies, and phase shifts into the Kuramoto equations which ultimately amount to affine changes of variables. For simplicity, we try to leave out these parameters.

Strogatz and Watanabe proved the existence of $N-3$ conserved quantities on the Kuramoto Model\cite{ws}. A priori, the Kuramoto Model is a vector field on the very high dimensional torus $(S^1)^N$, which is difficult to study directly. But by restricting to the submanifold where the $N-3$ quantities determined by the initial conditions are constant, we can reduce these equations to a three dimensional dynamical system. The conserved quantities are the cross ratios $\lambda_{ijkl} =  \frac{|e^{i\theta_i} - e^{i \theta_k}|}{|e^{i\theta_i} - e^{i \theta_l}|} \frac{|e^{i\theta_j} - e^{i \theta_k}|}{|e^{i\theta_j} - e^{i \theta_l}|}.$ There are ${N \choose 4}$ cross ratios, but it turns out only $N - 3$ are independent of each other. Observe that the fact that this ratio is conserved over time means that bodies which are incident at $t = 0$ remain incident for all time (this is trivial) but bodies which are not incident $t = 0$ will never collide because a nonzero cross ratio will always stay nonzero.

\section{Conformal Geometry Preliminaries}

A concise yet thorough textbook of Riemannian Geometry accessible to pure and applied mathematicians is John Lee's \textit{Riemannian Manifolds: An Introduction to Curvature}\cite{lee}.

A \textbf{Riemannian metric} on a smooth manifold $M$ is a section of the tensored cotangent bundle $(T^*)^{\otimes 2}M$ that is nondegenerate at every fiber. Informally, a Riemannian metric smoothly varies the inner product between tangent spaces. Two Riemannian manifolds are \textbf{conformal} if there is a diffeomorphism which scales each inner product by a positive factor. Conformal maps preserve angles and local orientation, but could possibly alter global geometric properties such as lengths of curves.

The most common nontrivial conformal map is \textbf{stereographic projection} of a sphere $S^d$ to $\R^{d+1}$. Embed $S^d$ in $\R^{d+1}$ in the standard way, which is via the subset $\{x \in \R^{d+1} \mid |x| = 1 \}$. Let $N = (0,0, \dots, 1) \in S^d$ be the north pole. For each $x \in S^d - N$, draw a line from $N$ to $x$ and map $x$ to the intersection of this line with the $x_{d+1} = 0$ plane, which we identify with $\R^d$. If we want to rigorously define the conformal map, we compactify $\R^d$ with one point and map $N$ to it. A critical observation is that the map acts as the identity on the equator, which is diffeomorphic to $S^{d-1}$. The stereographic projection also maps the southern hemisphere to the interior of the unit disc in $\R^{d}$ and the northern hemisphere to the exterior of the unit disc. Showing  sterographic projection is conformal is an elementary exercise in differential geometry, a proof of which is in Lee's book.

Since stereographic projection acts as the identity on the equator, the image of a Kuramoto model which takes place on the equator is invariant under this transformation, even with respect to distances. We will see that a generalized Kuramoto model on $S^d$ is described by a path of conformal maps on the ambient space $\R^{d+1}$, which we can then apply the inverse stereographic projection to get a path of conformal transformations on the sphere $S^{d+1}$.

In this paper, we concern ourselves with conformal maps of spheres, with the Riemannian metric inherited by the standard embedding. Rigid rotations of the sphere are conformal, but there are many others. Fractional linear transformations of the form $F(z) = \frac{az + b}{cz + d},$ with $ac-bd > 0$ are conformal maps of the Riemann sphere $\hat{\C} \cong S^2$. The formula is applied on the extended complex plane, but can be visualized as diffeomorphisms on the sphere via the inverse stereographic projection. The group of fractional linear transformations is called \textbf{the M\"obius Group}.

\section{Interpreting the Classical Results in Terms of Conformal Geometry}

We saw that the cross ratios between bodies on the Kuramoto model is preserved by its flow. A standard result in complex analysis states that the maps which preserve cross ratios are the M\"obius transformations. The space of the Kuramoto model is $S^1$, which we can embed in the standard way into $\C$, and its image under the inverse stereographic projection is the equator of $S^2$. 

Mirollo et. al. saw that the flow of a Kuramoto model can be described by a dynamical system on the Lie group of M\"obius transformations. However, this dynamical system did not lie in the entire Lie group, but rather the open subgroup of M\"obius transformations which preserve the unit disc and its boundary. If we apply the inverse stereographic projection, these are the conformal transformations of $S^2$ which preserve the southern hemisphere. These transformations take the form $M_{\zeta,w}(z) = \zeta ( \frac{w - z}{1 - \ov{w}z})$, where $|\zeta| = 1$ and $|w| < 1$. Hence, this group is diffeomorphic to the filled in torus $S^1 \times D^2$. This means if $\theta_i$ is a body on the Kuramoto model, its trajectory is modeled by $\theta_i = M_{\zeta,w}(\theta_i(0))$, where $\zeta$ and $w$ are parameters which are time dependent, but apply for all bodies. 

Consider the one dimensional Kuramoto Model defined by $A,B$ and $C$ as mentioned in the beginning. Let $a = -C+iB$. Then Mirollo et. al. wrote down the explicit equations on $S^1 \times D^2$ in the $(\zeta,w)$ coordinates in terms of our given data by differentiating the trajectory equation $\theta_i = M_{\zeta,w}(\theta_i(0))$ with respect to time. Specifically,

\begin{align*}
    \dot{w} &= -\frac{1}{2}(1 - |w|^2)\ov{\zeta}a(p) \\
    \dot{\zeta} &= iA(p)\zeta - \frac{1}{2}(\ov{w}a(p) - w \ov{a(p)}\zeta^2).
\end{align*}

This reduces an $N$ dynamical system on $(S^1)^N$ to just a three dimensional system, independent of $N$. This is a much more tractable system of equations!

\section{Generalized M\"obius Transformations}

The generalized Kuramoto model with $N$ bodies $x_1, x_2, \dots, x_N \in S^d$ takes the form $\dot{x_i} = \Omega x_i + X - \langle X, x_i \rangle x_i$, where $\Omega$ is a predetermined antisymmetric matrix independent of time and $i$, and $X \in \R^d$ is also independent of $i$, but possibly time dependent. Lohe saw the cross ratios $\lambda_{ijkl} = \frac{|x_i - x_k|}{|x_i - x_l} \frac{|x_j - x_k|}{|x_j - x_l|}$ are preserved. These quantities are preserved by M\"obius transformations on $\hat{\R}^{d+1} \cong S^{d+1}$ which map $S^d$ to itself. By deriving a formula for the M\"obius transformations of this form and then differentiating with respect to time, Lohe derived the reduced dynamical system. For completeness, we provide a derivation of the formula for these general M\"obius transformations, which is conspicuously absent in many texts on conformal geometry.

Liouville's Theorem from conformal geometry (which is unrelated to his theorem on bounded entire functions) states that for $d \geq 3$, conformal mappings on $\hat{\R}^{d+1}$ take the form $f(x) = b + \frac{\alpha A(x - a)}{|x - a|^\vep}$, where $a,b \in \R^{d+1}$, $\alpha \in \R$, $A \in SO_+(d+1)$, and $\vep = 0$ or $\vep = 2$. As mentioned before, the mappings which describe the Kuramoto model do not range over all possible conformal mappings, but rather are restricted to those which map $S^d$ to itself.

Rather than directly applying Liouville's Theorem, we can take the parametrization of the desired conformal maps in the $d = 2$ case and generalize. Recall that M\"obius transformation preserving the unit disc in $\R^2$ has the form $$M_{\zeta,w}(z) = \zeta \big( \frac{w - z}{1 - \ov{w}z} \big), |w| < 1$$

Note that some authors have $z - w$ in the numerator instead, which means the derivations here can have a similar, but different form as those which appear in other works. Navigating these conventions ultimately boils down to a transition map on coordinate charts describing the same manifold. Observe we have a Lorentz boost determined by the $w$ parameter, followed by a rotation of the $\zeta$ parameter. All conformal mappings can be decomposed as a product of this form. In higher dimensions, we replace multiplication by $\zeta$ with the rotation determined by a parameter in $SO(d)$, which we still call $\zeta$, which we know how to parametrize in theory (take the matrix exponentials of the vector space of antisymmetric $d \times d$ matrices). So we will focus on the Lorentz boost.

As with most complex fractions, we multiply the numerator and denominator by the denominator's complex conjugate.

\begin{align*}
    \frac{w - z}{1 - \ov{w}z} \cdot \frac{1 - w\ov{z}}{1 - w\ov{z}} & = \frac{w - w^2\ov{z} + w|z|^2}{|1 - \ov{w}{z}|^2}
\end{align*}

With the denominator, we can rewrite

\begin{align*}
    |1 - \ov{w}z|^2 &= \langle 1 - \ov{w}z, 1 - \ov{w}z \rangle \\
    &= \langle 1, 1 \rangle - 2\langle 1, \ov{w}{z} \rangle + \langle \ov{w}{z}, \ov{w}{z} \rangle
\end{align*}

Remember that we are taking the inner product of these complex numbers as vectors in $\R^2$, so $1$ actually denotes the vector $(1,0)$ and the inner product is the standard one on $\R^2$, not the Hermitian product because the components are real, not complex. Writing $w = w_1 + iw_2, z = z_1 + iz_2$, we can see
\begin{align*}
    \langle (1,0), \ov{w}{z} \rangle &= \tx{Re}(\ov{w}z)\\
    &=  \tx{Re}((w_1 - iw_2)(z_1 + iz_2))\\
    &= w_1z_1 + w_2z_2 \\
    &= \langle w, z \rangle.
\end{align*}

Thus,

\begin{align*}
    \langle (1,0), (1,0) \rangle - 2\langle 1, \ov{w}{z} \rangle + \langle \ov{w}{z}, \ov{w}{z} \rangle &= 1 - 2 \langle w, z \rangle + |\ov{w}z|^2 \\
    &= 1 - 2\langle w,z \rangle + |w|^2|z|^2 \\
    &= |w - z|^2 + (1 - |w|^2)(1 - |z|^2).
\end{align*}

For the denominator to be zero, we must have $w = z$ and $|w| = |z| = 1$, which cannot happen because $|w| < 1$ by assumption. Therefore, the denominator is never zero.

For the numerator, we make a clever addition and subtraction. Observe,

\begin{align*}
    w - w^2\ov{z} - z + w|z|^2 &= w - w^2\ov{z} - z + w|z|^2 - |w|^2(-z + w|z|^2) + |w|^2(-z + w|z|^2) \\
    &= (-z + w|z|^2)(1 - |w|^2) + w - w^2\ov{z} - z|w|^2 + w|w|^2|z|^2.
\end{align*}

The ultimate goal is to factor out a $w$ in the latter four terms, and the remainder just so happens to equal the denominator, leading to a nice cancellation. Notice that

\begin{align*}
    w - w^2\ov{z} - z|w|^2 + w|w|^2|z|^2 &= w - w^2\ov{z} - zw\ov{w} + w|w|^2|z|^2\\
    &= w(1 -w\ov{z} - z\ov{w} + |w|^2|z|^2).
\end{align*}

Further observe

\begin{align*}
    w\ov{z} + z\ov{w} &= \ov{\ov{w}z} + \ov{w}z\\
    &= 2\tx{Re}(\ov{w}z) \\
    &= 2\langle w, z \rangle.
\end{align*}

So cancellation yields a single $w$ summand, leading us to conclude the vector form of M\"obius transformations appearing in Lohe's paper:

$$M_{\zeta,w}(z) = \zeta\Big( \frac{(-z + w|z|^2)(1 - |w|^2)}{1 - 2\langle w,z \rangle + |w|^2|z|^2} + w \Big).$$

In the following, we ignore $\zeta$ and focus on the Lorentz boost. If we substitute $z' = \frac{-z}{|z|^2}$, we get $M_{w}(z) = \frac{(w + z')(1 - |w|^2)}{|w + z'|^2} + w $. This substitution simplifies the formula whilst allowing $M_w = I$ for $w = 0$. If $z = 0$, there is a removable discontinuity and we can see $M_{w}(0) = w$. We now show that $M_{\zeta,w}$ maps the unit ball to itself \cite{blue}.

Observe
\begin{align*}
    |M_w(z)|^2 &= \langle M_w(z), M_w(z) \rangle \\
    &= |w|^2 + 2\frac{(1-|w|^2)}{|w+z'|^2} \langle w, w + z' \rangle + \Big(\frac{1 - |w|^2}{|w + z'|^2} \Big)^2 |w+z'|^2.
\end{align*}

Clearing out denominators, we get

\begin{align*}
    |w-z'|^2 |M_w(z)|^2 &= |w|^2|w-z'|^2 + 2(1 - |w|^2)\langle w, w + z' \rangle + (1 - |w|^2)^2 \\
    &= |w|^2(|w|^2 + 2 \langle w, z' \rangle + |z'|^2) + 2(1 - |w|^2)(|w|^2 + \langle w, z' \rangle) + (1 - 2|w|^2 + |w|^4) \\
    &= 2 \langle w, z' \rangle + 1 + |w|^2|z|^2 \text{ (a ton of things cancel)}\\
    &= ( 2 \langle w, z' \rangle + |w|^2 + |z'|^2) + (1 + |w|^2|z'|^2 - |w|^2 - |z|^2) \\
    &= |w + z'|^2 + (1 - |w|^2)(1 - |z'|^2).
\end{align*}

Therefore, we have that $|M_w(z)|^2 - 1 = \frac{(1 - |w|^2)(1 - |z'|^2)}{|w + z'|^2}$. From the assumptions that $|z| \leq 1$ and $|w| < 1$, we have that the numerator of the righthand side is nonpositive. Therefore, $|M_w(z)| \leq 1$, with equality when $|z| = 1$. 

This generalization is well-studied in conformal geometry, and it satisfies many of the same properties of the classical M\"obius transformations in two dimensions. The proofs that these maps respect properties enjoyed by the M\"obius transformations ultimately comes down to more high school algebra. Proofs can be found in many sources, such as Stoll's text. From now on, we will use the variable $R$ instead of $\zeta$ to emphasize the fact that it is a rotation matrix.

The properties we use are:

\begin{itemize}
    \item The maps $M_{R,w}(z)$ are diffeomorphisms of the open unit ball and its boundary sphere in $\R^{d}$. 

    \item These maps preserve the sphere metric on the unit disc induced by its inclusion in to the southern hemisphere of the sphere. 
    
    \item These maps also preserve the cross ratios $\lambda_{ijkl} = \frac{|z_i - z_k|}{|z_i - z_l|}\frac{|z_j - z_k|}{|z_j - z_l|}$. (Here $z_i$ refers to a body on the sphere $S^d$, not a vector component). In the two dimensional case, division by vectors makes sense because of the complex structure, but in general we have to take ratios of lengths.
    
    \item The inverse of $M_{R,w}$ is $(M_{-w} \circ R^T) = M_{R^T, R w}$. Recall that $R$ is orthogonal, so $R^{-1} = R^T$. 
    
    \item All maps with the above properties can be realized in the above form, so $SO(d) \times D^{d}$ is indeed a parametrization of the group of conformal maps on $S^{d+1}$ which preserve $S^d$.
    
    \item $SO(d) \times D^{d}$ is isomorphic to the Lorentz group $SO_+(d,1)$.

\end{itemize}

\section{Lohe's Reduced Equations}

The generalization of the Kuramoto Model of a collection of coupled bodies $\{x_i\}$ on $S^d$ is $\dot{x_i} = \Omega x_i + X - \langle x_i, X \rangle x_i$, where $\Omega$ is an antisymmetric frequency matrix and $X$ determines how the bodies are coupled. Lohe has shown this system has conserved quantities $\lambda_{ijkl} = \frac{|x_i - x_k|}{|x_i - x_l|}\frac{|x_j - x_k|}{|x_j - x_l|}$, which know are preserved by generalized M\"obius transformations. 

Hence, the trajectory of a body $x_i$ is $x_i(t) = M_{\zeta(t),w(t)}(x_i(0))$, and we can apply the same M\"obius transformation to every body to determine the state of the system at a given time. This reduces our dynamical system on $(S^d)^N$ with dimension $Nd$ to a much simpler dynamical system on a Lie group of dimension $(d+1) + \frac{(d-1)d}{2}$, which is independent of $N$. We are free to reparametrize $SO(d)$ and $D^{d+1}$ to time-dependent coordinates which happen to induce even nicer looking conserved quantities.

Let $u_i(t) = M_{w(t)}(x_i(t))$. If we assume $\dot{w} = \Omega w + \frac{1}{2}( 1 + |w|^2)X - \langle w, X \rangle w$,  we have $\langle u_i(t), u_j(t) \rangle$ as a conserved quantity for all $i$ and $j$. This is the so-called ``Watanabe-Strogatz" transform which changes coordinates such that the bodies appear to be moving under standard rigid motions of the sphere. This allows us to consider the rotation parameter separately. Such a $w$ can be derived by assuming the inner product is preserved and working backwards to find the constraint on $w$. This equation is independent of $N$, so taking $N \to \infty$ gives Ed Ott's differential equation for probability distribution on spheres. Thus, we have that $u_i$ is modeled by a trajectory of the form $u_i(t) = R(t)u_i(0)$, where $R(t)$ is a path in $SO(d)$. Without loss of generality, we may assume $R(0) = I$ via a Lie group translation.

Observe that $\dot{u_i} = \dot{R}u_i(0) = \dot{R}R^Tu_i(t)$. Let $A = \dot{R}R^T - \Omega$. The tangent space of $SO(d)$ is the space of antisymmetric matrices, so $\dot{R}$ is antisymmetric. Differentiating the identity $RR^T = I$, we can see $\dot{R}R^T = -R^T\dot{R} = (\dot{R}R^T)^T$, so we can conclude $A$ is also antisymmetric. Now we have the equation differential $\dot{R} = (\Omega + A)R$, and Lohe shows that choosing $A_{ij} = X_iw_j - X_jw_i$ satisfies $\dot{R}R^Tu_i = (A + \Omega)Ru_i$, which can be verified via a computer algebra system. 

\section{Next Steps}

If we can solve for $R$ and $w$, we can then write the M\"obius transformations in the following time dependent matrix form:

\begin{align*}
    M_{R,w} &=
    \begin{bmatrix}
    R & f(\vec{w}) \\
    -f(\vec{w})^T & 1 \\
    \end{bmatrix}
    \\
    \end{align*}
    
where $f$ is a canonical transformation of $w$ to a vector in $\R^{d}$ with a well-known formula. Then if we want to track the movement of a body over time, we apply this time dependent matrix to its initial position.

Max Lohe derived the reduced equations, and we reinterpreted this result in terms of dynamical systems on conformal mapping groups. Since the stereographic projection acts as the identity on the equatorial sphere, we simultaneously have two topologically equivalent vector fields which we can visualize on $\R^{d+1}$ and $S^{d+1}$. The theory of vector fields on spheres comes with many topological obstructions which could lead to further understanding of the qualitative behavior of the Kuramoto model.

\end{document}